\documentclass[12pt,twoside]{amsart}

\usepackage{amsmath, amssymb, amscd, paralist,tabularx,supertabular,amsmath, verbatim, amsthm, amssymb, mathrsfs, manfnt,times,latexsym,amscd,graphicx}
\usepackage{pstricks}
\usepackage[all]{xy}

\usepackage{fullpage}

\DeclareMathOperator{\SL}{SL}

\newcommand{\Q}{\mathbb Q}
\newcommand{\R}{\mathbb R}
\newcommand{\Z}{\mathbb Z}

\newcommand{\bfF}{\mathbb F}

\newcommand{\calD}{\mathcal{D}}

\newcommand{\GG}{\mathrm G}
\newcommand{\HH}{\mathrm H}
\newcommand{\cO}{\mathcal{O}_k}

\numberwithin{equation}{section}

\theoremstyle{plain}

\newtheorem{theorem}{Theorem}[section]

\newtheorem{prop}[theorem]{Proposition}
\newtheorem{conj}[theorem]{Conjecture}
\newtheorem{prob}[theorem]{Problem}

\theoremstyle{remark}

\newtheorem*{ack}{Acknowledgements}
\title{Counting isospectral manifolds}

\author{Mikhail Belolipetsky}\thanks{Belolipetsky is partially supported by CNPq and FAPERJ research grants.}
\address{IMPA\\
Estrada Dona Castorina 110\\
22460-320 Rio de Janeiro, Brazil}
\email[] {mbel@impa.br}

\author{Benjamin Linowitz}\thanks{Linowitz is partially supported by an NSF Mathematical Sciences Postdoctoral Fellowship.}
\address{Department of Mathematics\\Oberlin College\\Oberlin, OH 44074\\USA}
\email[] {benjamin.linowitz@oberlin.edu}

\begin{document}

\begin{abstract}
Given a simple Lie group $H$ of real rank at least $2$ we show that the maximum cardinality of a set of isospectral non-isometric $H$-locally symmetric spaces of volume at most $x$ grows at least as fast as $x^{c\log x/ (\log\log x)^2}$ where $c = c(H)$ is a positive constant.
In contrast with the real rank $1$ case, this bound comes surprisingly close to the total number of such spaces as estimated in a previous work of Belolipetsky and Lubotzky \cite{BL}. Our proof uses Sunada's method, results of \cite{BL}, and some deep results from number theory. We also discuss an open number-theoretical problem which would imply an even faster growth estimate.
\end{abstract}

\maketitle

\section{Introduction}

Two compact Riemannian manifolds are called isospectral if the associated Laplace-Beltrami operators have the same eigenvalue spectrum. The notion of isospectrality extends naturally to the setting of compact locally symmetric Riemannian orbifolds (see \cite{Gor}). Recall that these orbifolds are obtained as the quotients of the symmetric space $X$ of a Lie group $H$ by cocompact discrete subgroups $\Gamma < H$. Two orbifolds $\Gamma_1\backslash X$ and $\Gamma_2\backslash X$ are isometric if and only if the lattices $\Gamma_1$ and $\Gamma_2$ are conjugate in the isometry group of $X$. We will denote by $\mathrm{IL}_{H}(x)$ the maximum cardinality of a set of conjugacy classes of pairwise isospectral cocompact irreducible lattices in $H$ whose covolume with respect to a fixed Haar measure is bounded above by $x$.

The question of the existence of isospectral non-isometric Riemannian manifolds was popularized in the nineteen-sixties by Marc Kac \cite{Kac} who famously asked ``Can one hear the shape of a drum?''. Although Kac was primarily interested in isospectral planar domains, Milnor \cite{Milnor} had already constructed isospectral non-isometric flat tori in dimension $16$. The first examples of isospectral non-isometric Riemann surfaces were constructed by Vign\'eras \cite{Vign} in 1980. A few years later Sunada \cite{Sunada} developed a group theoretic approach to the problem which led to many new examples. In \cite{BrGG}, Brooks, Gornet and Gustafson applied Sunada's method so as to construct arbitrarily large sets of pairwise isospectral non-isometric Riemann surfaces and then studied the growth of the size of these sets as a function of volume. In the notation above, their result says that for $H = \SL(2, \R)$ we have $\mathrm{IL}_{H}(x) \ge x^{c\log x}$ for some constant $c>0$. It is worth noting that this lower bound grows faster than any polynomial function of $x$. In contrast to this, it was shown in \cite{Lin} that the number of different isospectral Riemann surfaces which can be constructed by Vign\'eras' method is bounded above by a polynomial function of volume. In \cite{McR}, McReynolds extended the lower bound of \cite{BrGG} to the groups of isometries of higher dimensional real hyperbolic spaces and complex hyperbolic $2$-space. He also showed that for any non-compact simple Lie group $H$ and every positive integer $r$, there exists a strictly increasing sequence $(x_j)$ such that $\mathrm{IL}_{H}(x_j) \ge x_j^r$.

In the present paper we obtain explicit super-polynomial lower bounds for the function $\mathrm{IL}_{H}(x)$ for almost all non-compact semisimple Lie groups $H$. Our main result is as follows.

\begin{theorem}\label{main_thm}
Let $H$ be a non-compact semisimple Lie group which has no simple factors of type $\mathrm{A}_1$ and which contains irreducible lattices. Then there exists a positive constant $c = c(H)$ such that
$$\mathrm{IL}_{H}(x) \ge x^{c\log x/ (\log\log x)^2},$$
and a bound of the same shape holds for torsion-free lattices.
\end{theorem}

In particular, Theorem~\ref{main_thm} holds for any simple Lie group of rank at least $2$. Our method unfortunately falls short of proving an analogous result for the remaining higher rank semisimple Lie groups containing irreducible lattices:
\begin{conj}
The lower bound of Theorem~\ref{main_thm} holds for semisimple Lie groups of type $\mathrm{A}_1$ which have more than one simple factor.
\end{conj}

Let us note that the groups of type $\mathrm{A}_1$ often appear more difficult for the questions about distribution of lattices. For example, these groups were excluded from the main theorem in \cite{BP} and were brought into its scope only much later in \cite{BGLS}.

The principal case for our theorem is when the group $H$ has real rank at least 2. To put it into a more general perspective, we may compare the growth rate of the number $\mathrm{IL}_{H}(x)$ of isospectral lattices with the total number $\mathrm{L}_{H}(x)$ of conjugacy classes of irreducible lattices in $H$ whose covolume is bounded above by $x$. It is known that for the group $H$ of isometries of hyperbolic $n$-space the function $\mathrm{L}_{H}(x)$ grows at least as fast as $x^{c'x}$ and that the same lower bound applies if we restrict ourselves to counting only cocompact arithmetic lattices \cite{BGLS}. Only a tiny proportion of this quantity corresponds to the pairwise isospectral lattices of Brooks--Gornet--Gustafson and McReynolds. In the case of higher rank Lie groups the situation is markedly different. By Margulis's arithmeticity theorem all irreducible lattices in higher rank semisimple Lie groups are arithmetic (cf. \cite{WitMor15}). The growth of $\mathrm{L}_{H}(x)$ for these groups was investigated in \cite{BL}, where it was shown that it is at rate $x^{a\log x}$, which is \emph{surprisingly close} to the lower bound from Theorem~\ref{main_thm}. It is quite possible that our lower bound for $\mathrm{IL}_{H}(x)$ is still not optimal and the actual growth reaches the rate $x^{c\log x}$ (the same type as the total growth of lattices!). This possibility leads to an intricate number-theoretical problem which we will discuss in the end of Section~\ref{cyclotomic seq}.

The proof of the main result of this paper is based on group theory and number theory. The theorem then applies to higher rank locally symmetric manifolds and tells us that many of them are mutually isospectral. The geometry of this phenomenon remains mysterious and could provide an interesting subject for future study.

\section{Sunada's Method}\label{sunada}
In this section we review Sunada's construction \cite{Sunada} of isospectral manifolds together with some generalizations and applications.

Given a finite group $G$ and element $g\in G$, denote by $[g]_{G}$ the $G$-conjugacy class of $g$. Two subgroups $H_{1}$ and $H_{2}$ of $G$ are said to be {\em almost conjugate} if
\[ \#\left(H_{1}\cap [g]_{G}\right)= \#\left(H_{2}\cap [g]_{G}\right)\]
for all $g\in G$.

Conjugate subgroups of $G$ are always almost conjugate, but we will be mainly interested in nonconjugate almost conjugate subgroups. This is motivated by the following theorem of Sunada:

\begin{theorem}(Sunada \cite{Sunada})\label{theorem:sunadastheorem}
Let $M$ be a closed Riemannian manifold, $G$ a finite group, and \[\varphi: \pi_{1}(M)\longrightarrow G\] a surjective homomorphism. If $H_{1}$ and $H_{2}$ are almost conjugate subgroups of $G$ then the manifolds associated to $\varphi^{-1}(H_{1})$ and $\varphi^{-1}(H_{2})$ are isospectral.
\end{theorem}

In order to construct different isospectral covers of a closed Riemannian manifold $M$, the theorem requires a surjective homomorphism $\pi_{1}(M)\rightarrow G$ onto a finite group $G$ containing almost conjugate nonconjugate subgroups. The proof of this paper's main result will make use of the following variant of Sunada's theorem:

\begin{theorem}\label{theorem:sunadastheoremvariant}
Let $M$ be a closed Riemannian orbifold, $G$ a finite group with almost conjugate subgroups $\{H_{i}\}_{i\in I}$, and $G'$ a finite group containing $G$. If \[\varphi: \pi_{1}(M)^{orb}\longrightarrow G'\] is a surjective homomorphism then the orbifolds associated to $\varphi^{-1}(H_{i})$ are pairwise isospectral covers of the orbifold associated to $\varphi^{-1}(G)$.
\end{theorem}

Here the change from the group $G$ to the larger group $G'$ is a minor variation which can be obtained by the same argument as \cite{Sunada} and the extension of the result to the orbifold setting follows from the work of B\'erard \cite{Ber1, Ber2}.

Among the smallest finite groups which admit almost conjugate nonconjugate subgroups are the semi-direct product $(\Z/8\Z) \rtimes (\Z/8\Z)^*$ and the group $\SL(3, \Z/2\Z)$. These groups were used to construct some nice explicit examples of isospectral manifolds. We refer the reader to \cite{Bus86} and \cite{BrTse} for the first such constructions. The groups $G$ for which we will employ Theorem \ref{theorem:sunadastheoremvariant} are Heisenberg groups over finite fields. These groups were first applied to the construction of isospectral manifolds in a paper of Brooks, Gornet and Gustafson \cite{BrGG}.

% In Section \ref{heisenberg} we will exhibit Heisenberg groups over finite fields which have large numbers of nonconjugate almost conjugate subgroups. This will allow us to employ Corollary \ref{theorem:sunadastheoremvariant} so as to produce families of pairwise isospectral manifolds. Proving that these manifolds are pairwise non-isometric is more nuanced and will require an argument of Belolipetsky and Lubotzky \cite{BL}. Finally, in Section~\ref{count_all} we will elaborate on another result from \cite{BL}, which would allow us to produce large families of non-isometric and isospectral lattices associated to the infinite Hilbert class field towers.

% In order to construct large families of isospectral nonisometric manifolds we will apply Corollary~\ref{theorem:sunadastheoremvariant} to products of Heisenberg groups over finite fields.

Let $q=p^n$ be a prime power and recall that the {\em Heisenberg group} $\HH(\mathbb F_q)$ is the subgroup of $\SL_3(\mathbb F_q)$ consisting of upper-triangular unipotent matrices; that is,
\[ \HH(\mathbb F_q)=\left\{ \begin{pmatrix}
1 & x & y \\
0 & 1 & z \\
0 & 0 & 1
\end{pmatrix}: x,y,z\in\mathbb F_q \right\}.\]

Consider the subgroup $H_1$ of $\HH(\mathbb F_q)$ given by
\[ H_1=\left\{ \begin{pmatrix}
1 & x & 0 \\
0 & 1 & 0 \\
0 & 0 & 1
\end{pmatrix}: x\in\mathbb F_q \right\}.\] It is clear that the order of $H_1$ is $q$, hence the order of any subgroup of $\HH(\mathbb F_q)$ which is almost conjugate to $H_1$ is also equal to $q$.

In \cite[p. 317]{BrGG} Brooks, Gornet and Gustafson exhibited a family of $p^{n(n-1)}$ almost conjugate nonconjugate subgroups of $\HH(\mathbb F_q)$, where $q=p^n$. These subgroups were constructed as the orbit of $H_1$ under a family of {\em almost linear automorphisms} of $\HH(\mathbb F_q)$.

Recall that the conjugacy classes in a direct product of groups are simply the product of conjugacy classes in the constituent groups. In particular if $\{G_i\}_{i\in I}$ is a family of finite groups and $\{H_{i,1},H_{i,2} < G_i\}_{i\in I}$ are almost conjugate nonconjugate subgroups of the $G_i$ then $\prod_{i\in I} H_{i,1}$ and $\prod_{i\in I} H_{i,2}$ are almost conjugate nonconjugate subgroups of $\prod_{i\in I} G_i$.

From this discussion we are now able to conclude the following.

\begin{prop}\label{prop31}
Let $q_1=p_1^{n_1},\dots, q_r=p_r^{n_r}$ be a sequence of prime powers and set \[G=\HH(\mathbb F_{q_1})\times\cdots\times \HH(\mathbb F_{q_r}).\] The group $G$ contains at least $p_1^{n_1(n_1-1)}\cdots p_r^{n_r(n_r-1)}$ almost conjugate nonconjugate subgroups of order $q_1\cdots q_r$.
\end{prop}

\section{Constructing isospectral lattices over a field}\label{count_k}

Let $H$ be a connected semisimple Lie group without compact factors and $\Gamma$ be an irreducible \emph{arithmetic subgroup} of $H$ defined over a number field $k$. Recall that this means there exists a simply connected, absolutely simple algebraic $k$-group $\GG$ which admits an epimorphism
\linebreak
${\phi:\GG(k\otimes_\Q\R) \to H}$ whose kernel is compact and such that $\phi(\GG(\cO))$ is commensurable with $\Gamma$, where $\cO$ denotes the ring of integers of $k$. We refer to the book of Witte Morris \cite{WitMor15} for more information about arithmetic subgroups and their properties.

We shall assume that $\Gamma = \phi(\Lambda)$ with $\Lambda \subset \GG(\cO)$. Let us fix a rational prime $p$ with decomposition $p\cO = \mathcal{P}_1^{e_1}\mathcal{P}_2^{e_2}\cdots \mathcal{P}_m^{e_m}$, $\cO/\mathcal{P}_i = \bfF_{p^{n_i}}$, for which the reduction map $\pi_p: \GG(\cO) \to \GG(\cO/p\cO)$ is surjective on $\Lambda$. Thus we have
\begin{equation}\label{eq41}
 \pi_p(\Lambda) = \GG(\cO/p\cO) = \GG(\bfF_{p^{n_1}}) \times \GG(\bfF_{p^{n_2}}) \times \ldots \times \GG(\bfF_{p^{n_m}}).
\end{equation}
The strong approximation property implies that this is true for all but finitely many primes $p$ (cf. \cite[Chapter~7]{PlR}). 
% In fact, what we shall need is a weaker property requiring that $\phi(\Lambda)$ contains a certain product of Heisenberg groups.

From now on let us assume that the group $\GG$ is not of the type $\mathrm{A}_1$. We would like to show that the group $G_0 = \GG(\bfF_{p^{n_1}}) \times \ldots \times \GG(\bfF_{p^{n_m}})$ contains a product of Heisenberg groups to which we can then apply Proposition~\ref{prop31}. This requires some standard theory of finite groups of Lie type for which we will refer to Carter's book \cite{Ca}.

Let $G = \GG(\bfF_{q})$ be a finite quasi simple group of Lie type of rank at least $2$, and let us assume that the characteristic of the field $\bfF_q$ is bigger than $3$. This assumption will allow us to exclude some special cases (Suzuki, Ree, and Tits groups) and does not in any way restrict our applications. Thus $\GG$ is a Chevalley group or a twisted group of type $^2\mathrm{A}_n$, $^2\mathrm{D}_n$, $^2\mathrm{E}_6$ or $^3\mathrm{D}_4$ (cf. \cite[Chapters 4 and 13]{Ca}). Let us take two positive roots $a$, $b$ such that $a+b$ is also a root and in the twisted cases the roots $a$, $b$ do not belong to the same equivalence class (see \cite[Sec. 13.2.1]{Ca} for the precise definition of the equivalence classes induced by twisting). It follows that the linear combinations of $a$, $b$ which are roots form a root subsystem $\Sigma$ of type $\mathrm{A}_2$, $\mathrm{B}_2$ or $\mathrm{G}_2$ (cf. \cite[Lemma 3.6.3]{Ca}).

Consider a unipotent subgroup $U < G$ generated by the following root subgroups:
\begin{itemize}
 \item[(1)] $X_a$, $X_b$, if  $\Sigma$ is of type $\mathrm{A}_2$;
 \item[(2)] $X_a$, $X_{a+b}$, if  $\Sigma$ is of type $\mathrm{B}_2$;
 \item[(3)] $X_b$, $X_{3a+b}$, if  $\Sigma$ is of type $\mathrm{G}_2$.
\end{itemize}
Using Chevalley's commutator formula it is easy to check that $U$ is isomorphic to a Heisenberg group defined over $\bfF_q$.

Now let $M$ be a compact locally symmetric space associated to a cocompact arithmetic lattice $\Gamma= \phi(\Lambda)$, i.e. $M = \Gamma\backslash H/K$ where $K$ is a maximal compact subgroup of $H$. Our considerations, combined with Proposition~\ref{prop31} and Sunada's theorem (Theorem \ref{theorem:sunadastheoremvariant}), show that $M$ has at least $p^{n_1(n_1-1)+ \cdots + n_m(n_m-1)}$ pairwise isospectral covers corresponding to nonconjugate subgroups of $\Gamma = \pi_1(M)^{orb}$ of index bounded above by the order of the group $G_0$ divided by $p^{n_1}\cdots p^{n_m}$.

The following proposition summarizes the discussion above.

\begin{prop}\label{prop41}
Let $\GG$ be an absolutely simple, simply connected algebraic group defined over a number field $k$ whose type is different from $\mathrm{A}_1$. Assume that for a rational prime $p >3$ the prime ideals in the decomposition of $p\cO$ in $\cO$ have inertia degrees $n_1$, $n_2$, \ldots, $n_m$, and that the group $\Lambda \subset \GG(\cO)$ maps onto $\GG(\cO/p\cO)$ under the reduction map. Then $\Lambda$ contains at least $p^{n_1(n_1-1)+ \cdots + n_m(n_m-1)}$ subgroups which are mapped to almost conjugate nonconjugate subgroups of $\GG(\cO/p\cO)$ of index bounded above by $\#\GG(\cO/p\cO)/(p^{n_1}\cdots p^{n_m})$.
\end{prop}

Note that although the subgroups provided by the proposition are mutually nonconjugate in $\Lambda$ they may be still be conjugate in the larger group $H$, which contains $\phi(\Lambda)$ as a lattice. Up to a constant factor, the sizes of $H$-conjugacy classes are equal to the conjugacy classes of lattices with respect to the action of the isometry group of $X$ (see \cite[Section~5]{BL} for more details).  Hence in order to count non-isometric isospectral manifolds we must control these conjugations. The desired control can be obtained by employing the results of \cite[Section~5]{BL}, assuming that we can construct sufficiently many almost conjugate subgroups. The latter means that we need at least one of the inertia degrees $n_i$ to be close to the degree of the extension $k/\Q$; ideally the prime $p$ should be inert in $k/\Q$. The main challenge is how to achieve this while keeping covolume of $\Gamma$ reasonably small. This is the goal of the next section.

\section{A number-theoretic construction} \label{cyclotomic seq}

In order to apply Proposition~\ref{prop41} and Sunada's method for constructing large families of isospectral lattices we would like to have a sequence $(k_i)_{i = 1, 2,\ldots}$ of extensions of a fixed number field $k$ which satisfy the following properties:
\begin{itemize}
\item[(1)] All of the fields $k_i$ have a given, fixed number of complex places;
\item[(2)] The degrees $d_i = [k_i:\Q]$ tend to infinity;
\item[(3)] There exists a prime ideal $\mathcal{P}$ in $\cO$ which stays \emph{inert} in all $k_i$.
\end{itemize}
As in \cite{BL}, if we would like the volume of the associated quotient spaces to grow slowly we should require:
\begin{itemize}
\item[(4)] The root discriminant $rd_i = \calD_{k_i}^{1/d_i}$ of $k_i$ grows slowly with $i$, where $\calD_{k_i}$ denotes the absolute value of the discriminant of the field $k_i$.
\end{itemize}

Ideally, in condition~(4) we would like to assume that the root discriminant is uniformly bounded. This assumption is achieved, for instance, by unramified field extensions arising from an infinite class field tower. The existence of these towers was established in a seminal paper of Golod and Shafarevich \cite{GS}, and the associated field extensions were used in \cite{BL} to prove an asymptotically sharp lower bound for the number of higher rank lattices of bounded covolume. The problem we encounter now is different because of condition~(3), which is hard to achieve in unramified extensions. Indeed, if we consider an infinite, full Hilbert class field tower (as in \cite{BL}) then it is easy to see that condition~(3) can never be achieved! This is, of course, not the only possible construction of infinite sequences of unramified extensions but the observation above illustrates the kind of difficulties that arise.  We will make a few more comments about this problem at the end of the section, by now we turn our attention to a different construction based on cyclotomic extensions.

We refer to \cite{Wash}   %\cite[Chaper~IV]{Lang}
for the basic theory of cyclotomic fields. Let us consider a sequence of cyclotomic fields of prime conductor $\Q(\zeta_{p_i})$, where $\zeta_{p_i}$ denotes a primitive $p_i$th root of unity and the $p_i$ run through the prime numbers greater than $2$.  In order to satisfy condition~(1) (with no complex places) we will consider the maximal real subfields $k_i = \Q(\zeta_{p_i} + \zeta_{p_i}^{-1})$ of the $\Q(\zeta_{p_i})$. We can then produce fields with any fixed number of complex places and similar asymptotic properties to the $(k_i)$ by taking appropriate quadratic extensions defined by Pisot numbers, as was done in \cite{BL}. We will come back to this construction in the next section.

Having taken $k_i = \Q(\zeta_{p_i} + \zeta_{p_i}^{-1})$, we now have
\begin{align*}
& d_i = [k_i:\Q] = \frac{p_i-1}2;\\
& \calD_{k_i} = p_i^{\frac12(p_i-3)}.
\end{align*}
Thus, the root discriminant $rd_i = p_i^{\frac{p_i-3}{p_i-1}}$ grows $\sim2d_i$, which is reasonably slow.

By \cite[Theorem~2.13]{Wash} a prime $p > 2$ is inert in $\Q(\zeta_{p_i})/\Q$ (and hence in $k_i/\Q$) if it is a primitive root modulo $p_i$. Artin's famous conjecture on primitive roots states that a given integer which is neither a perfect square nor $-1$ is a primitive root modulo infinitely many primes. Although this conjecture is still open, there are some positive, unconditional results in its direction which will be sufficient for our purposes. For example, Heath-Brown proved that there are at most two exceptional primes for which Artin's conjecture fails \cite{HB}. Although Heath-Brown's theorem does not give any extra information about the exceptional primes, it implies that one or more of $5$, $7$, and $11$ is a primitive root modulo $\ell$ for infinitely many primes $\ell$. By restricting our sequence of conductors $p_i$ to these primes $\ell$ we obtain a sequence of field extensions $k_i$ of the field $k = \Q$ which satisfy the conditions (1)--(4). We summarize the main properties of this construction for future reference:

\begin{prop}\label{prop51}
There exists an infinite sequence of totally real fields $k_i$ of degree $d_i \to \infty$ with root discriminant $rd_i = (2d_i+1)^{1-\frac1{d_i}}$ in which a prime $3 < p \le 11$ stays inert.
\end{prop}

Let us come back to the question about a possible improvement upon this result. This suggests the following open problem:
\begin{prob}\label{prob1}
Do there exist infinite sequences of number fields $k_i$ which satisfy conditions (1)--(3) and have uniformly bounded root discriminant?
\end{prob}

A positive solution to Problem \ref{prob1} would allow us to show that $\mathrm{IL}_H(x) \ge x^{c\log x}$ by applying the argument of Section~\ref{count_all}. Unfortunately, as was mentioned above, condition~(3) and uniform boundedness in condition~(4) are difficult to achieve simultaneously. If the $k_i$ constitute an infinite (full) Hilbert class field tower, i.e. every $k_i$ is the Hilbert class field of $k_{i-1}$, then all of the extensions $k_i/k$ are Galois. Were a prime ideal $\mathcal{P}_i$ of $k_i$ to be inert in two successive levels $k_{i+1}$ and $k_{i+2}$ of the tower, the corresponding Galois group $\mathrm{Gal}(k_{i+2}/k_i)$ would be the cyclic group generated by the Frobenius element $\left[\frac{k_{i+2}/k_i}{\mathcal{P}_i}\right]$, hence abelian, contradicting the fact that these are two distinct levels in the class field tower. However, there is still a possibility here to build upon infinite class field towers. We recall that by a result of Hajir \cite{Haj}, the $p$-ranks of the class groups in an infinite $p$-class field tower grow to infinity. This implies that there are many possible choices of infinite sequences of $p$-extensions of $k$ that are contained in the maximal unramified $p$-extension $k^{(\infty)}$. The finiteness argument above does not apply here any more and it is possible that the required construction exists. The problem reduces to a subtle analysis of inertia and capitulation of ideals in unramified extensions. We conclude with a remark that a positive solution to Problem~\ref{prob1} may shed new light on the structure of the pro-$p$ Galois group $\mathrm{Gal}(k^{(\infty)}/k)$ which plays an important role in number theory.

\section{Counting isospectral lattices}\label{count_all}

The method that we are going to use in this section is based on the proof of the lower bound for the number of lattices $\mathrm{L}_\HH(x)$ given in \cite[Sections~6 and 8]{BL}. We will review the construction here, skipping some of the technical parts for which we will refer the reader to \cite{BL}.

We begin with a sequence of number fields $k_i$ and extensions $l_i$ of the $k_i$ that satisfy
\begin{equation}\label{asm_k_l}
\calD_{k_i} \le (c_0d_i)^{d_i},\quad  \calD_{l_i}/\calD_{k_i}^{[l_i:k_i]} \le (c'_0d_i)^{d_i},
\end{equation}
where $\calD$ denotes the absolute value of the field discriminant, $d_i = [k_i:\Q]$ is the degree of the field $k_i$, and $c_0$, $c'_0$ are positive constants. Moreover, we require that the fields $k_i$ and $l_i$ have their signature determined by $H$. For example, if $H$ is a real Lie group of type $\mathrm{B}_n$, then we have to assume that all the fields $k_i$ are totally real and $l_i = k_i$. The other groups may require for every $i$ the field $l_i$ to be a quadratic extension of $k_i$ and both of them having a certain number of complex places (see \cite[pp.~3136, 3144]{BL}). In the basic totally real case we take the fields $k_i$ provided by Proposition~\ref{prop51}. The general case can be reduced to this by a clever trick using Pisot numbers (see \cite[Lemma~3.4]{BL}, \cite[Corollary~3.5]{BL} and the pages cited above). In all the cases the base fields have degree at most $2$ over the fields from Proposition~\ref{prop51}.

The assumptions about the signature of the fields allow us to associate for every pair $(k, l) = (k_i, l_i)$ an absolutely simple simply connected $k$-group $\GG$ such that:
\begin{itemize}
\item[(1)] $\GG(k\otimes_\Q\R)$ admits an epimorphism to $H$ whose kernel is compact;
\item[(2)] $\GG$ is quasi-split over $k_v$ for every $v\in V_f\smallsetminus \{ v_0 \}$, where $V_f$ denotes the set of finite places of $k$;
\item[(3)] The quasi-split inner form of $\GG$ splits over $l$.
\end{itemize}
This is done using Prasad--Rapinchuk's extension of the Borel--Harder theorem in \cite{PR06} (see also \cite[Proposition 6.1]{BL}). The place $v_0$ in condition (2) can be an arbitrary selected finite place of $k$, thus we can assume that its residue characteristic is greater than $11$.

The next step is to construct arithmetic subgroups of $\GG(k)$ whose covolume and reduction mod $p\cO$ we can control. To this end, for every finite place $v$ different from $v_0$ we can choose a ``nice'' parahoric subgroup $\mathrm{P}_v$ (which is always a special parahoric and is hyperspecial whenever $l$ is unramified over $k$ at $v$). Together with an arbitrarily chosen parahoric subgroup $\mathrm{P}_{v_0}$, we obtain a coherent collection $\mathrm{P} = (\mathrm{P}_v)_{v \in V_f}$ of parahoric subgroups of $\GG$. We now take $\Lambda = \GG(k)\cap\prod_{v\in V_f} \mathrm{P}_v$, the corresponding principal arithmetic subgroup of $\GG(k)$. The reader who is not familiar with the theory of the algebraic groups over non-archimedean local fields can think of this construction as a generalization of the basic example $\SL_n(\Z) = \SL_n(\Q)\cap \prod_{\text{primes }p}\SL_n(\Z_p)$. Now the covolume of a principal arithmetic subgroup $\Lambda$ can be evaluated using Prasad's volume formula \cite{Pr}. Together with the assumptions \eqref{asm_k_l} on the fields $k$ and $l$ and some elementary estimates, this gives us an upper bound for the covolume of $\Lambda' = \phi(\Lambda) < H$ (cf. \cite[pp.~3137--3138]{BL}):

\begin{align*}
\mu(H/\Lambda') &\le (c_1d^{\frac12(\mathrm{dim}(\GG)+s)})^{d},\\
\text{where } c_1 &= c_0^{\frac12\mathrm{dim}(\GG)}{c'_0}^{\frac12s}
\prod_{i=1}^{r}\frac{m_i!}{(2\pi)^{m_i+1}}\: p_0^{\mathrm{dim}(\GG)}
\left(\frac{\pi^2}{6}\right)^{r}.
\end{align*}
To review the notation, recall that $\mathrm{dim}(\GG)$, $r$ and $m_i$ denote the dimension, rank and Lie exponents of the group $\GG$, $s$ is an integer constant also determined by $\GG$, $p_0$ is the characteristic of the reduction field at $v_0$, and the other constants are as before.

Finally note that by the construction the reduction map $\pi_p: \Lambda \to \GG(\cO/p\cO)$ is surjective for any rational prime $p\neq p_0$.

We now choose $p' \neq p_0$ to be the inert prime provided by Proposition~\ref{prop51}. The inertia degree of $p'$ in the field $k$ is greater than or equal to $\frac12d$ (as $k$ has degree at most $2$ over the corresponding $k_i$ in the proposition). By Proposition~\ref{prop41}, the group $\Lambda$ has at least $p'^{\frac{d}2(\frac{d}2-1)}$ subgroups associated to almost conjugate nonconjugate subgroups of index at most $c_3^d$, where $c_3 = p'^{\mathrm{dim}(\GG) - \frac12}$ (we used the fact that $\#\GG(\cO/p\cO) \le p'^{d\cdot\mathrm{dim}(\GG)}$). This gives $c_4^{d^2}$ isospectral lattices in $H$ of covolume at most $(c_1c_3d^\gamma)^d$, with $\gamma = \frac12(\mathrm{dim}(\GG)+s)$ (where we can take $c_4 = p'^{\frac14-\epsilon}$ for a small $\epsilon > 0$ assuming that $d \ge \frac{1}{2\epsilon}$).

We still have to check that sufficiently many of these lattices are nonconjugate in $H$. To this end we recall the results of \cite[Section~5]{BL}. All subgroups constructed above contain the $p'$-congruence subgroup $\Lambda(p')$ of $\Lambda$. By \cite[Corollary~5.3]{BL}, each conjugacy class contains at most $c_3^d (c_1d^\gamma)^{d\cdot C} p'^{d\cdot\mathrm{dim}(\GG)} \le (c_5d^{\gamma\cdot C})^d$ such subgroups. Hence we have at least $c_3^{d^2}(c_5d^{\gamma\cdot C})^{-d}$ different conjugacy classes of lattices of covolume at most $(c_1c_3d^\gamma)^d$.

We now deduce the bound in Theorem~\ref{main_thm}. Denote the upper bound for covolume by $x$. We have:
\begin{align*}
& \log x = d\log(c_1c_3d^\gamma) \ge \gamma d\log d; \\
& \log\log x \ge \log d.
\end{align*}
For the logarithm of the number of nonconjugate lattices, assuming $d$ is large enough, we have:
\begin{align*}
& \log(c_3^{d^2}(c_5d^{\gamma\cdot C})^{-d}) = d^2\log c_3 - d \log(c_5d^{\gamma\cdot C}) \ge c_6 d^2,
\end{align*}
hence the number itself is $\ge c_7^{d^2}$ and  we can write it as
\begin{align*}
& c_7^{(\log x)^2/(\log(c_1c_3d^\gamma))^2} \ge x^{a\log x/(\log\log x)^2},
\end{align*}
assuming $x$ is large enough.

This finishes the proof of the theorem for lattices that may contain elements of finite order. By Selberg's lemma, any lattice in a linear group contains a torsion free congruence subgroup of finite index. In \cite{BE}, Belolipetsky and Emery suggested a variant of this result which allows one to obtain a good upper bound for the index of the torsion-free subgroup. The idea of the construction in \cite{BE} is that instead of a principal congruence subgroup corresponding to a sufficiently large prime, we can choose an intersection of preimages of $p$-Sylow subgroups of the quotients at two different primes (see \cite[Section~2]{BE} for the details). This means that we can replace the arithmetic lattices $\Lambda$ by torsion-free subgroups of index at most $c^d$ whose local structures differ only at the places which lie over two fixed rational primes, say, $p = 2$ and $3$. With these lattices in hand we can repeat the above argument for the lower bound of the number of isospectral non-isometric $H$-manifolds. \qed

\begin{ack}
We are grateful to Farshid Hajir, Jonah Leshin, Ben McReynolds and Ilir Snopche for helpful discussions. We also thank the anonymous referee for a careful proofreading of the paper. 
\end{ack}


\begin{thebibliography}{11}

\bibitem{BGLS} M. Belolipetsky, T. Gelander, A. Lubotzky, A. Shalev, Counting arithmetic lattices and surfaces,
{\em Ann. of Math. (2)} {\bf 172} (2010), 2197--2221.

\bibitem{BE} M. Belolipetsky, V. Emery, Hyperbolic manifolds of small volume, {\em Doc. Math.}, {\bf 19} (2014), 801--814.

\bibitem{BL} M. Belolipetsky, A. Lubotzky,
\newblock \emph{Manifolds counting and class field towers},
\newblock Adv. Math. \textbf{229} (2012), no. 6, 3123-3146.

\bibitem{Ber1} P. B\'erard, Transplantation et isospectralit\'e. I,  {\em Math. Ann.}, {\bf 292} (1992), 547--559.

\bibitem{Ber2} P. B\'erard, Transplantation et isospectralit\'e. II, {\em J. London Math. Soc. (2)}, {\bf 48} (1993), 565--576.

\bibitem{BP} A. Borel, G. Prasad, Finiteness theorems for discrete subgroups of bounded covolume
in semi-simple groups, {\em Inst. Hautes \'Etudes Sci. Publ. Math.} {\bf 69}
(1989), 119--171. Addendum: {\em ibid.,} {\bf 71} (1990), 173--177.

\bibitem{BrGG} R.~Brooks, R.~Gornet, and W.~H. Gustafson.
\newblock Mutually isospectral {R}iemann surfaces.
\newblock {\em Adv. Math.}, \textbf{138} (1998), no. 2, 306--322.

\bibitem{BrTse} R. Brooks, R. Tse, Isospectral surfaces of small genus, {\em Nagoya Math. J.}, {\bf 107} (1987), 13--24.

\bibitem{Bus86} P.~Buser, Isospectral Riemann surfaces, {\em Ann. Inst. Fourier}, {\bf 36} (1986), 167--192.

\bibitem{Ca} R. W. Carter,
\newblock \emph{Simple Groups of Lie Type},
\newblock Pure and Applied Mathematics 28, Wiley, London, 1972.

\bibitem{GS} E. S. Golod, I. R. Shafarevich, % I. R. \v Safarevi\v c,
On the class field tower, {\em Izv. Akad. Nauk SSSR Ser. Mat.}, {\bf 28}
(1964), 261--272 [Russian].

\bibitem{Gor} C. Gordon. Orbifolds and their spectra,
in {\em  Spectral geometry}, 49--71.
{\em Proc. Sympos. Pure Math.}, {\bf 84}, Amer. Math. Soc., Providence, RI, 2012.

\bibitem{Haj} F. Hajir, On the growth of $p$-class groups in $p$-class field towers, {\em J. Algebra}, {\bf 188} (1997), 256--271.

\bibitem{HB} D. Heath-Brown, Artin's conjecture for primitive roots, {\em Quart. J. Math. Oxford Ser.}, {\bf 37} (1986), 27--38.

\bibitem{Kac} M. Kac, Can one hear the shape of a drum? {\em Amer. Math. Monthly}, {\bf 73} (1966), 1--23.

% \bibitem{Lang} S. Lang, {\em Algebraic number theory}, Addison-Wesley, 1970.

\bibitem{Lin} B. Linowitz, Families of mutually isospectral Riemannian orbifolds, {\em Bull. Lond. Math. Soc.}, {\bf 47} (2015), 47--54.

\bibitem{McR} D. B. McReynolds, Isospectral locally symmetric manifolds, {\em Indiana Univ. Math. J. (2)}, {\bf 63} (2014), 533--549.

\bibitem{Milnor} J.~Milnor, 
\newblock \emph{Eigenvalues of the {L}aplace operator on certain manifolds}, 
\newblock Proc.\ Nat.\ Acad.\ Sci.\ U.S.A. \textbf{51} (1964), 542.

\bibitem{PlR} V. P. Platonov, A. S. Rapinchuk, {\em Algebraic groups and number theory}, Pure Appl. Math. {\bf 139}, Academic Press, Boston, 1994.

\bibitem{Pr} G. Prasad, Volumes of $S$-arithmetic quotients of semi-simple groups,
{\em Inst. Hautes \'Etudes Sci. Publ. Math.,} {\bf 69} (1989), 91--117.

\bibitem{PR06} G. Prasad, A. S. Rapinchuk, On the existence of isotropic forms
of semi-simple algebraic groups over number fields with prescribed local
behavior, {\em Adv. Math.}, \textbf{207} (2006), 646--660.

\bibitem{Sunada} T.~Sunada,
\newblock \emph{Riemannian coverings and isospectral manifolds},
\newblock Ann. of Math. (2), \textbf{121} (1985), 169--186.

\bibitem{Vign} M. F. Vign\'eras, Vari\'et\'es riemannienes isospectrales et non isometriques, {\em Ann. Math.}, {\bf 112} (1980), 21--32.

\bibitem{Wash} L.~Washington, {\em Introduction to Cyclotomic Fields}, Graduate Texts in Mathematics, vol. 83, Springer,
New York, 1997.

\bibitem{WitMor15} D.~{Witte Morris}.
\newblock \emph{Introduction to Arithmetic Groups}.
\newblock Deductive Press, 2015.

\end{thebibliography}
\end{document}